# Limiting behavior of solutions of multidimensional Landau-Lifshitz equations with second approximation of effective field (I)[*]

Ganshan Yang[1,*]   Boling Guo[2]   Jihui Wu[3]

1) Institute of Mathematics, Yunnan Normal University, Kunming P. O. 650092, China.

Department of Mathematics, Yunnan Nationalities University, Kunming P. O. 650031, China.

2) Institute of Applied Physics and Computational Mathematics, Nonlinear Center for Studies

of China Academy of Engineering, P. O. Box 8009, Beijing 100088, China

3) Department of Mathematics, Huainan Normal Univerisity, Huainan P. O. 232038, China

Corresponding author: Ganshan Yang, email: ganshanyang@yahoo.com.cn

**Abstract**

By a maximum principle under global control map from $R^{n+1}$ to $S^2$, $\delta$-viscosity supersolution and subsolution of multidimensional Landau-Lifshitz equations(pert I: hot harmonic case) with second approximation of effective field are built. Utilizing the $\delta$-viscosity supersolution and subsolution, the viscosity solution of the equations is showed and the limiting behavior of the solution of the equations is obtained. Exactly, there exist two disjoint open subsets such that the solutions tend to the point $(0,1,0)$ on arbitrary compact subsets in one of them and tend to the point $(0,-1,0)$ on arbitrary compact subsets in the other of them, respectively.

**MR(2000) subject Classification** 15A15, 15A09, 35Q55

## 1 Introduction

In general, the Landau-Lifshitz equations is given by

$$u_t = u \times H^e,\ |u| = 1,\ n = 3. \quad (1.a)$$

This equation describes motion of the magnetization $u$ around the effective field $H^e$, which is derived from some energy functional $E = E(u)$ according to

$$H^e = -\frac{\delta E}{\delta u}.$$

[*]The work is supported by National Natural Science Foundation of China (No. 11161057), the Institute of Mathematical Sciences of Chinese University of Hongkong, Institute of Mathematics of Fudan University and Institute of Mathematics of Peking University.

Therefore the effective field itself is a function of the magnetization and equation (1.a) is, in general, highly nonlinear. An explicit form of the Landau-Lifshitz equations is obtained by specifying the energy functional $E$. The total magnetic energy has the form:

$$E = E_{an}(u) + E_{ex}(u) + E_{ext}, \quad (1.b)$$

where $E_{an}(u)$ is an anisotropy energy, $E_{ex}(u)$ is the exchange energy and $E_{ext}$ is the energy from the external magnetic field $H$. The equations proposed by Landau and Lifshitz in 1935 (see [18]) may be written as the following form:

$$\begin{cases} \frac{\partial u}{\partial t} = \lambda_1 u \times H^e - \lambda_2 u \times (u \times H^e) & in \quad \Omega \times (0,T), \\ H^e := -H(u) + \sum_{i,j=1}^{3} \frac{\partial}{\partial x_j}(a_{ij} \frac{\partial u}{\partial x_i}) + H & in \quad \Omega \times (0,T), \end{cases} \quad (1.c)$$

where $\lambda_1$ $\lambda_2$ are constants, $\lambda_2 > 0$.

In 1985, A. Visintin in [26] studied the energy functional $E$ which corresponds to (1.c), he showed that the corresponding minimization problem has at least one solution in the following case:

$$\begin{cases} H(u) \text{ is gradient of a convex function } R^3 \to R^+, \\ \{a_{ij}\} \text{ is a postive definite } 3 \times 3 \text{ matrix.} \end{cases}$$

Whether the solutions of minimization problem on the energy function are smooth? There are some results on one dimensional case in [38]. On the multidimensional case, the existence of weak solution was first proved in 1986 by Sulem P.-L., Sulem C., Bardos C. (see [25]) and Zhou, Guo (see [37]), respectively. Sulem P.-L., Sulem C., Bardos C. proved that for any $S_0$ with $|S_0(x)| = 1$ a.e. and $\partial S/\partial x_i \in L^2(R^d)$, there exists a global weak solution of $\partial S/\partial t = S \times \Delta S$, $S(x,0) = S_0(x)$. Zhou and Guo proved the global existence of weak solution for generalized Landau-Lifhsitz equations without Gilbert term. In 2000 Nai-Heng Chang, J. Shatah and K. Uhlenbeck considered initial value problem for the 2-dimensional cylindrical symmetric Landau-Lifshitz equations (see [4]). There are some papers to consider the partial regularity and the existence of weak solutions (see[7, 20, 22, 29]), to analysis the concentration sets (see[21]), to study the relation between weak solution and harmonic map (see[13, 14]), to discuss the stability(see[36, 17, 16, 30]). Recent some years, there are many important programs on Schrödinger map equation by dealing more subtly with derivative term (how deal with the derivative Schrödinger, we suggest to refer [27, 28]), global solution(see [2]), the equivalence under the Hasimoto transformation(see [1]) and blow up solution(see [23]).

Whether there exists global smooth solution for (1.c) is an open problem. Indeed, the multidimensional Landau-Lifshitz equation is hard to solve exactly, even if two dimensional Heisenberg's equations. Only when there is a lack of an external field or an anisotropy field, many authors had attempted Hirota's famous method, Lie algebraic structure, and Bäcklund transform, Double transform, etc., but exact solution were still not constructed.

In 2000-2001 Guo, Han, Yang constructed several exact solutions to the cylindrical symmetric $n$-dimension Landau-Lifshitz equation in [11, 12]. In 2001, Guo, Yang constructed a family of exact disc solutions to the two-dimension Landau-Lifshitz equation in [15].

In 2003, Yang and Chang obtained an exact transform between with and without an external field(see [32]). By the transform, some sufficient and necessary conditions that the solution with an external field converges to a certain solution without an external field are given(see [32, 35]).

When both external field and anisotropy field are present there is not any paper discussing the explicit solution of the multidimensional Landau-Lifshitz equation until 2008. In 2009, Yang and Guo researched the multidimensional Landau-Lifshitz equation with an uprush external field and anisotropy field, and gave a method to find explicit solution, provided some explicit solutions for two dimensional case (see[33]). Mu and Liu et. gave a method to find an external field and anisotropy field such that the solution remains the same(see [24]).

In 2011, Yang and Liu present a method to solve the Landau-Lifshitz equations with the external magnetic field and anisotropy field and a new family of symmetric solutions, i.e. spherical cone symmetric solution, show how these solutions evolve as time, give a series of explicit dynamic spherical cone symmetric solutions. They find an interesting phenomenon that a equiv-variant solution is static (see [34]). In 2012, [31] provided an explicit blow up solution of Schrödinger equation derived from Schrödinger map, showed the non-equivalence between the Schrödinger equation and Landau-Lifshitz equation and also found that two class of equivariant solutions of Landau-Lifshitz equation are static.

In order to find out a way of studying (1.c), to doing research on some especial cases are very useful. In a neighborhood of Curie temperature, as a first approximation it is assumed that

$$H(u) = \frac{\partial \Phi(u)}{\partial u},$$

where

$$\Phi(u) = \sum_{i,j=1}^{3} b_{ij} u_i u_j$$

with symmetric and positive definite tensor $\{b_{ij}\}$. Therefore it is important to study the existence of solution of the following three-dimension Landau-Lifshitz equations with first approximation of anisotropy field $H(u)$:

$$\begin{cases} \frac{\partial u}{\partial t} = -u \times (u \times (\triangle u + H - \lambda H(u))) + \gamma u \times (\triangle u + H - \lambda H(u)) \\ \qquad\qquad \text{in } \Omega \times (0, \infty), \\ u = g \quad \text{on } \partial\Omega \times (0, \infty), \\ u \in S^2 \text{ in } \Omega \times (0, \infty), \end{cases} \quad (1.d)$$

where $\lambda$ is a positive parameter, $\gamma$ is a nonnegative parameter, $g = (g_1, g_2, g_3) \in C^{2,\alpha_0}(\partial\Omega; R^3) \cap (\partial\Omega; S^2)$ $(0 < \alpha_0 < 1)$.

The equation (1.d) is a form of Landau-Lifshitz equations presented in 1935 (see [18]). It describes the evolvement of spin field $u$ in continuum ferromagnets. $\triangle u + H - \lambda H(u)$ is the effective field containing contributions from interaction of fields. $H$ is the demagnetize field and $H(u)$ is the anisotropy field induced by $u$. If there are at least two components of $H(u)$ to be identically vanished, $H(u)$ is called as uni-direct effective field. Otherwise $H(u)$ is called as multi-direct effective field. The case with uni-direct effective field $H(u) = (0, 0, 2u_3)$ and boundary condition $g = (g_1, g_2, 0)$ was studied by Zhai Jian in [36]. He gave the existence of stable static solution of (1.d). We studied the case with first approximation of multi-direct effective field (see [16]).

However, the first approximation is not very satisfactory (see [19, 26]). Thus it is necessary to study the higher order approximation of anisotropy field $H(u)$. In this paper we discuss the second order approximation of the anisotropy field $H(u) = (\alpha u_1 u_3^2, \alpha u_2 u_3^2, \alpha u_1^2 u_3 + \alpha u_2^2 u_3 + \beta u_3)$ that which agrees with the one proposed in [19, 26]. We discover that if the boundary condition still take the form $g = (g_1, g_2, 0)$, then the solution of (1.d) with above the second order approximation of the anisotropy field still is stable; but if the boundary condition take the form as $g = (g_1, g_2, g_3)(g_3 \neq 0)$, then the stability is not uniform with respect to $\lambda$, even the case of $\gamma = 0$

$$\begin{cases} \frac{\partial u}{\partial t} = -u \times (u \times (\triangle u - \frac{1}{\varepsilon^2} H(u))) & in \quad R^n \times (0, T), \\ u \in S^2 & in \quad R^n \times (0, T), \end{cases} \quad (1.1)$$

where $\varepsilon > 0$ is a small parameter, $H(u) = (\alpha u_1 u_3^2, \alpha u_2 u_3^2, \alpha u_1^2 u_3 + \alpha u_2^2 u_3 + \beta u_3)$, $\alpha + \beta > 0$, $\beta > 0$, $\alpha, \beta$ are both real constants. $\Omega$ be a bounded domain in $R^n$ with $C^3$ boundary, $S^2$ denotes the unit sphere in $R^3$, $u = (u_1, u_2, u_3) : \Omega \to S^2$ is a three-dimensional vector valued unknown function and "$\times$" denotes the vector cross product in $R^3$.

Notice that the harmonic heat flow problem is seen as a special case of (1.1) without perturbation. Next Proposition 2.1 enlightens us to consider in another view the relationship between (1.d) and the mean curvature motion as some authors have considered Allen-Cahn equation by using viscosity solutions in [3, 5, 6, 8, 9, 10]. The notion of viscosity solution for nonlinear first-order partial differential equation was introduced by M. G. Crandall and P. L. Lions in [6], this notion has been proved quite successful in the existence and uniqueness theory of solutions of Hamilton-Jacobi equations(see [6]) and the singular Limit of Ginzhurg-Landau equation(see [3]). We expected such notion to be applied to Landau-Lifshitz equation. However the method of viscosity solution is valid only for signal equation, such as the Allen-Cahn equation is a singular function and not three dimensional vector, but Landau-Lifshitz equations is combined of three equations. How build viscosity solution for a system? It is very vagal to put directly the common viscosity solution of all equations in a system as a viscosity solution of this system (even supersolution or subsolution). Thus we built an "approximate" viscosity supersolution and subsolution: $\delta$-viscosity supersolution and subsolution for one of Landau-Lifshitz type equations with second approximation of effective field. Consequently we study the singular limit of the vector-valued Landau-Lifshitz equation by using those relationships between mean curvature motion and the $\delta$-viscosity supersolution and subsolution.

For the $n$ dimensional Landau-Lifshitz equations with three-dimensional vector values in $u \in S^2$, we did not know its limiting behavior. Our difficulty is mostly that the $n$ dimensional Landau-Lifshitz equations with three-dimensional vector values in $u \in S^2$ is a system, and as the viscosity solutions of every equation in a system may not be common viscosity solutions of all equations. Therefore it is very difficult to study the limiting behavior of $n$ dimensional Landau-Lifshitz equations by viscosity solutions. To overcome such difficulty we first give some concepts, e.g. $\delta$-viscosity super-solution and subsolution (see Definition 2.2-2.4). Then we build a maximum principle by a global control map $v^\varepsilon(r(x,t))$ from $R^{n+1}$ to $S^2$.

The present paper consists of eight parts. In section 2, we define the $\delta$-viscosity solutions,

supersolution and subsolution of (1.1).

In section 3 we first provide some further properties of signed distance function (see Lemma 3.1 and 3.2), then construct an auxiliary function (see(3.7)) and give some its properties (see Lemma 3.1-Remark 3.6).

In section 4, in order to obtain the $\delta$-viscosity supersolution and subsolution of (1.1), we construct a global control map $v^\varepsilon(r(x,t))$ from $R^{n+1}$ to $S^2$ by a smooth solution of an ordinary differential equations relative to (1.1). It is the foundation of this paper. In view of our purpose, we have to find special solutions with some special properties for the ordinary differential equations (see (4.3-4.4)). It is the most difficult to find out the solutions satisfying our purpose for the ordinary differential equations. Therefore, there are complicated calculations and some techniques in this section.

In section 5, we construct a pair of explicit $\delta$-viscosity supersolution and subsolution for (1.1) by the auxiliary function and those Lemmas obtained in section 3, this is one of our prime results.

In section 6, we provide a maximum principle and a comparison principle under the global control map, and prove the existence and uniqueness of $\delta$-viscosity solution. This is also one of our prime results.

In section 7, using the characteristic of $\delta$-viscosity supersolution and subsolution and those results built in the section 2-6, we obtain our prime result on the limiting behavior of the solution of the $n$ dimensional Landau-Lifshitz equations with three-dimensional vector value in a unit sphere, and prove that there exist two disjoint open subsets such that the solutions of the multidimensional Landau-Lifshitz equations with values in a unit sphere tend to the point $(0, 1, 0)$ and the point $(0, -1, 0)$ on arbitrary inner compact sets in one and another of the two open subsets, respectively.

In section 8, the method of $\delta$-viscosity supersolution and subsolution enlightens us to obtain a directly method of viscosity supersolution and subsolution of multidimensional Landau-Lifshitz equations with the second approximation of effective field.

## 2  Definition of $\delta$-viscosity solutions

Next, we use the notations $u_\xi$ and $u_{j\xi}$ to denote the partial derivative of the three-dimensional vector $u$ and the partial derivative of jth component in three-dimensional vector $u$, respectively, etc. To introduce some conceptions called as $\delta$-viscosity solutions of (1.1), we first prove the following proposition.

**Proposition 2.1** *Let $u \in C^{2,\alpha_0;\, 1,\frac{\alpha_0}{2}}(R^n \times [0,T]; R^3)$, then $u$ is a solution of (1.1) if and only if $u$ is also a solution of the following equations*

$$\begin{cases} u_2 u_{1t} - u_1 u_{2t} - (u_2 \triangle u_1 - u_1 \triangle u_2) = 0 \\ u_{3t} - \triangle u_3 - |\nabla u|^2 u_3 + \frac{1}{\varepsilon^2} u_3(1-u_3^2)(\alpha + \beta - 2\alpha u_3^2) = 0 \\ u \in S^2 \text{ in } R^n \times (0,T) \end{cases} \quad (2.1)$$

*where*

$$|\nabla u|^2 = |\nabla u_1|^2 + |\nabla u_2|^2 + |\nabla u_3|^2,$$

*provided that*
$$|u_3| \neq 1, \ for \ all \ (x,t) \in R^n \times (0,T).$$

Notice that the condition $|u_3| \neq 1, \ for \ all \ (x,t) \in R^n \times (0,T)$. It insures (2.1) to be equivalent to (1.1) in the meaning of $C^\infty(R^n \times (0,T))$.

We seek some solutions with the form $u = v(r(x,t))$ for (1.1), where $v(r) \in C^2(R; R^3) \cap S^2$, $r(x,t) \in LSC(R^n \times [0,\infty))$. Here $LSC(A)$ and $USC(A)$ denote the bounded lower and upper semicontinuous functions on $A$, respectively. When $r(x,t) \in C^2(R^n \times [0,\infty))$, $u$ is a smooth solution of (1.1).

Assume that $r(x,t) \in C^2(R^n \times [0,\infty))$, then

$$u_{it} = v_{ir}r_t, \quad \Delta u_i = v_{irr}|\nabla r|^2 + v_{ir}\Delta r, \quad |\nabla u|^2 = |v_r|^2|\nabla r|^2. \tag{2.2}$$

Substitute (2.2) into (2.1), we can get

$$\begin{cases} (v_2 v_{1r} - v_1 v_{2r})(r_t - \triangle r) - (v_2 v_{1rr} - v_1 v_{2rr})|\nabla r|^2 = 0 \\ v_{3r}(r_t - \triangle r) - (v_{3rr} + |v_r|^2 v_3)|\nabla r|^2 + \frac{1}{\varepsilon^2}v_3(1-v_3^2)(\alpha+\beta-2\alpha v_3^2) = 0 \\ v \in S^2 \end{cases} \tag{2.3}$$

The equation above enlightens us to introduce the concepts of $\delta$-viscosity solution of (1.1) as below.

Now, let $0 < \delta \ll 1$ be fixed.

**Definition 2.2** *A three dimensional vector valued function $u$ is called a $\delta$-viscosity supersolution of (1.1), provided that there is a three dimensional vector valued function (global control map) $v \in C^2(R; R^3) \cap S^2$ and a function $r \in LSC(R^n \times [0,\infty))$ such that $u(x,t) = v(r(x,t))$, for all $(x,t) \in R^n \times [0,\infty)$; and make the $v(r)$ be a viscosity supersolution of the first equation in (2.3) for $r \notin (\frac{\delta}{4}, \frac{\delta}{2})$ and the second equation in (2.3) for $r \in R$.*

The definitions on viscosity solution (or supersolution, subsolution) of singed equation were given in [5] and [9].

**Definition 2.3** *A three dimensional vector valued function $u$ is called a $\delta$-viscosity subsolution of (1.1), provided that there is a three dimensional vector valued function (global control map) $v \in C^2(R; R^3) \cap S^2$ and a function $r \in USC(R^n \times [0,\infty))$ such that $u(x,t) = v(r(x,t))$, for all $(x,t) \in R^n \times [0,\infty)$; and make the $v(r)$ be a viscosity subsolution of the first equation in (2.3) for $r \notin (-\frac{\delta}{2}, -\frac{\delta}{4})$ and the second equation in (2.3) for $r \in R$.*

**Definition 2.4** *A three dimensional vector valued function $u$ is called a $\delta$-viscosity solution of (1.1), if it is both $\delta$-viscosity supersolution and $\delta$-viscosity subsolution of (1.1). Specially $0$-viscosity solution (supersolution, subsolution, respectively) is called a viscosity solution (supersolution, subsolution, respectively).*

**Remark 2.5** *Let $v(s) \in C^2(R; R^3) \cap S^2$ and $|v_3| \neq 1$ for all $s \in R$, then $v(r)$ is a viscosity solution (supersolution, subsolution, respectively) of (1.1) iff $r$ is a common viscosity solution (supersolution, subsolution, respectively) of two equations of (2.3).*

**Remark 2.6** *From Definition 2.2-2.4 if $u = v(r(x,t))$ is a δ-viscosity supersolution (or solution, subsolution, respectively) of (1.1), then $v(r(x,t))$ is a δ-viscosity supersolution (or solution, subsolution, respectively) of (2.3). Roughly the δ-viscosity supersolution (or solution, subsolution, respectively) of a system is the viscosity supersolution (or solution, subsolution, respectively) of an equation in this system and the approximate viscosity supersolution (or solution, subsolution, respectively) of the other equations in this system.*

In fact, our δ-viscosity supersolution implies the following:

**Remark 2.7** *Let three dimensional vector valued function $u$ is a δ-viscosity supersolution of (1.1), then for each pair $\phi \in C^2(R^n \times [0, \infty))$ and $(x_0, t_0) \in R^n \times (0, \infty)$, if $r - \phi$ has a local minimum 0 at the point $(x_0, t_0)$, we have the following results at the point $(x_0, t_0)$:*

(A)
$$v_2(\phi)v_1''(\phi) - v_1(\phi)v_2''(\phi) = 0, \ when \ v_2(\phi)v_1'(\phi) - v_1(\phi)v_2'(\phi) = 0; \qquad (2.4)$$

$$\phi_t - \triangle \phi \geq \frac{v_2(\phi)v_1''(\phi) - v_1(\phi)v_2''(\phi)}{v_2(\phi)v_1'(\phi) - v_1(\phi)v_2'(\phi)} |\nabla \phi|^2, \ if \ only \ \phi \notin (\frac{\delta}{4}, \frac{\delta}{2}), \qquad (2.5)$$
$$when \ v_2(\phi)v_1'(\phi) - v_1(\phi)v_2'(\phi) \neq 0.$$

(B)
$$(v_3''(\phi) + |v'(\phi)|^2 v_3(\phi))|\nabla \phi|^2 - \tfrac{1}{\varepsilon^2} v_3(\phi)(1 - v_3^2(\phi))(\alpha + \beta - 2\alpha v_3^2(\phi)) = 0, \ when \ v_3'(\phi) = 0;$$

$$\phi_t - \triangle \phi \geq \frac{v_3''(\phi) + |v'(\phi)|^2 v_3(\phi)}{v_3'(\phi)} |\nabla \phi|^2 - \frac{\alpha + \beta - 2\alpha v_3^2(\phi)}{\varepsilon^2} \frac{v_3(\phi)(1 - v_3^2(\phi))}{v_3'(\phi)}, \ when \ v_3'(\phi) \neq 0. \qquad (2.6)$$

From Remark 2.7, we see that if

$$v_3'(\phi) \neq 0 \ and \ v_2(\phi)v_1'(\phi) - v_1(\phi)v_2'(\phi) \neq 0$$

at the point $(x_0, t_0)$, we have the following results at the point $(x_0, t_0)$:

$$\phi_t - \triangle \phi \geq \max \left\{ \frac{v_2(\phi)v_1''(\phi) - v_1(\phi)v_2''(\phi)}{v_2(\phi)v_1'(\phi) - v_1(\phi)v_2'(\phi)} |\nabla \phi|^2, \right.$$
$$\left. \frac{v_3''(\phi) + |v'(\phi)|^2 v_3(\phi)}{v_3'(\phi)} |\nabla \phi|^2 - \frac{\alpha + \beta - 2\alpha v_3^2(\phi)}{\varepsilon^2} \frac{v_3(\phi)(1 - v_3^2(\phi))}{v_3'(\phi)} \right\},$$

for $\phi \notin (\frac{\delta}{4}, \frac{\delta}{2})$;

$$\phi_t - \triangle \phi \geq \frac{v_3''(\phi) + |v'(\phi)|^2 v_3(\phi)}{v_3'(\phi)} |\nabla \phi|^2 - \frac{\alpha + \beta - 2\alpha v_3^2(\phi)}{\varepsilon^2} \frac{v_3(\phi)(1 - v_3^2(\phi))}{v_3'(\phi)},$$

for $\phi \in (\frac{\delta}{4}, \frac{\delta}{2})$.

δ-viscosity subsolution implies the following:

**Remark 2.8** *Let three dimensional vector valued function $u$ is a δ-viscosity subsolution of (1.1), for each pair $\phi \in C^2(R^n \times [0, \infty))$ and $(x_0, t_0) \in R^n \times (0, \infty)$, if $r - \phi$ has a local maximum 0 at the point $(x_0, t_0)$, we have the following results at the point $(x_0, t_0)$:*

(A)
$$v_2(\phi)v_1''(\phi) - v_1(\phi)v_2''(\phi) = 0, \ when \ v_2(\phi)v_1'(\phi) - v_1(\phi)v_2'(\phi) = 0;$$

$$\phi_t - \triangle \phi \leq \frac{v_2(\phi)v_1''(\phi) - v_1(\phi)v_2''(\phi)}{v_2(\phi)v_1'(\phi) - v_1(\phi)v_2'(\phi)} |\nabla \phi|^2, \ if\ only\ \phi \notin (-\frac{\delta}{2}, -\frac{\delta}{4}), \quad (2.7)$$
$$when\ v_2(\phi)v_1'(\phi) - v_1(\phi)v_2'(\phi) \neq 0.$$

(B)

$$(v_3''(\phi) + |v'(\phi)|^2 v_3(\phi))|\nabla \phi|^2 - \frac{1}{\varepsilon^2} v_3(\phi)(1 - v_3^2(\phi))(\alpha + \beta - 2\alpha v_3^2(\phi)) = 0, \ when\ v_3'(\phi) = 0;$$

$$\phi_t - \triangle \phi \leq \frac{v_3''(\phi) + |v'(\phi)|^2 v_3(\phi)}{v_3'(\phi)} |\nabla \phi|^2 - \frac{\alpha + \beta - 2\alpha v_3^2(\phi)}{\varepsilon^2} \frac{v_3(\phi)(1 - v_3^2(\phi))}{v_3'(\phi)}, \ when\ v_3'(\phi) \neq 0.$$

From Remark 2.8, we see that if

$$v_3'(\phi) \neq 0\ and\ v_2(\phi)v_1'(\phi) - v_1(\phi)v_2'(\phi) \neq 0$$

at the point $(x_0, t_0)$, we have the following results at the point $(x_0, t_0)$:

$$\phi_t - \triangle \phi \leq \min \left\{ \frac{v_2(\phi)v_1''(\phi) - v_1(\phi)v_2''(\phi)}{v_2(\phi)v_1'(\phi) - v_1(\phi)v_2'(\phi)} |\nabla \phi|^2, \right.$$
$$\left. \frac{v_3''(\phi) + |v'(\phi)|^2 v_3(\phi)}{v_3'(\phi)} |\nabla \phi|^2 - \frac{\alpha + \beta - 2\alpha v_3^2(\phi)}{\varepsilon^2} \frac{v_3(\phi)(1 - v_3^2(\phi))}{v_3'(\phi)} \right\},$$

for $\phi \notin (-\frac{\delta}{2}, -\frac{\delta}{4})$;

$$\phi_t - \triangle \phi \leq \frac{v_3''(\phi) + |v'(\phi)|^2 v_3(\phi)}{v_3'(\phi)} |\nabla \phi|^2 - \frac{\alpha + \beta - 2\alpha v_3^2(\phi)}{\varepsilon^2} \frac{v_3(\phi)(1 - v_3^2(\phi))}{v_3'(\phi)},$$

for $\phi \in (-\frac{\delta}{2}, -\frac{\delta}{4})$.

Notice that (2.5) may not hold for $\delta$-viscosity supersolution when $\phi \in (\frac{\delta}{4}, \frac{\delta}{2})$. Similarly, (2.7) may not hold for $\delta$-viscosity subsolution when $\phi \in (-\frac{\delta}{2}, -\frac{\delta}{4})$. This is essential difference from viscosity supersolution. Its signification will be seen in the sections 4-8.

## 3 Auxiliary functions

For the sake of completeness, we also recall the definition and two lemmas concerning signed distance function, but the lemmas here are slightly different from original ones. As $\partial \Omega$ is smooth, we may choose a continuous function $h$ to be smooth with $|\nabla h| = 1$ near $\partial \Omega$, and $h$ is a constant outside some ball so that

$$h(x) \begin{cases} > 0 & if\ x \in \Omega, \\ = 0 & if\ x \in \partial \Omega, \\ < 0 & if\ x \in R^n \setminus \overline{\Omega}. \end{cases} \quad (3.1)$$

Assume that $\partial \Omega$ is compact, then according to Theorem 3.2 and 4.2 in [9], we see that there is an unique weak solution $w(x, t) \in C(R^n \times [0, \infty))$ of mean curvature partial differential equation (this is a function, not three dimensional vector ) with initial value $h$

$$\begin{cases} w_t = (\delta_{ij} - \frac{w_{x_i} w_{x_j}}{|\nabla w|^2}) w_{x_i x_j} & in\ R^n \times (0, \infty), \\ w = h & in\ R^n \times \{t = 0\}. \end{cases} \quad (3.2)$$

Now define

$$I_t \equiv \{x \in R^n \mid w(x, t) > 0\},$$

$$\Gamma_t \equiv \{x \in R^n \mid w(x,t) = 0\},$$
$$O_t \equiv \{x \in R^n \mid w(x,t) < 0\},$$
$$I \equiv \{(x,t) \in R^n \times (0,\infty) \mid w(x,t) > 0\},$$
$$O \equiv \{(x,t) \in R^n \times (0,\infty) \mid w(x,t) < 0\}.$$

According to Theorem 7.6 in [9], there exists a time $t^* > 0$ such that $\Gamma_t$ is the boundary of a nonempty open set for $0 \leq t \leq t^*$. Let us assume that $t^* = \inf\{t > 0 \mid \Gamma_t = \emptyset\}$,

$$d(x,t) = \begin{cases} dist(x, \Gamma_t) & if \ x \in I_t, \\ 0 & if \ x \in \Gamma_t, \\ -dist(x, \Gamma_t) & if \ x \in O_t, \end{cases} \quad (3.3)$$

for $x \in R^n$, $0 \leq t \leq t^*$, where $d$ is called a signed distance function (see [8]).

**Lemma 3.1** *For every $t \in [0, t^*)$, $\Gamma_t$ is a nonempty compact set, $d$ is Lipschitz continuous in $x$ and $|\nabla d| = 1$ for almost all $x \in R^n$.*

**Lemma 3.2** *Let $\phi \in C^2(R^n \times [0, \infty))$ and $(x_0, t_0) \in R^n \times (0, t^*)$ satisfy that*

$$d - \phi \ \text{has a local minimum at the point} \ (x_0, t_0), \quad (3.4)$$

*then*

$$|\nabla \phi| = 1, \quad \text{at the point} \ (x_0, t_0).$$

To guarantee that (2.3) has $\delta$-viscosity supersolution, we construct a smooth auxiliary function $\eta : R \to R$ satisfying the following conditions (in fact, these conditions are necessary for the proof of Theorem 7.1):

(i) $\eta \in C^2(R)$, $0 \leq \eta'(z) \leq 1$ and $0 \leq \eta''(z) \leq \frac{C}{\delta}$.

(ii) There exist positive constants $\tau_1$, $\tau_2$ and $a$ that depending only on $\delta$ and $t^*$ and satisfying $\frac{\delta}{4} < \tau_1 < \tau_2 < \frac{\delta}{2}$ and $0 < a < 1$ such that $\eta \in C^3(\tau_1, \frac{\delta}{2})$,

$$\begin{array}{lll} \eta'(z) = 0, & for & z \leq \frac{\delta}{4}; \\ 0 < \eta'(z) \leq a < 1, & for & z \in (\frac{\delta}{4}, \tau_2); \\ \eta'(z) = 1, & for & z \geq \frac{\delta}{2}; \\ \eta'''(z) < 0, & for & z \in (\tau_1, \frac{\delta}{2}). \end{array} \quad (3.5)$$

For example, we assume that

$$\eta(z) = \begin{cases} C_1(\delta), & z \leq \frac{\delta}{4}, \\ \text{the polynomial with order four}, & \frac{\delta}{4} \leq z \leq \frac{\delta}{2}, \\ z + C_2(\delta), & z \geq \frac{\delta}{2}. \end{cases} \quad (3.6)$$

Applying the method of undetermined coefficients, we can prove the existence of the $\eta$ satisfying above conditions. For instance $\eta$ may be taken the following

$$\eta(z) = \begin{cases} -\frac{5}{8}\delta, & z \leq \frac{\delta}{4}, \\ -\frac{32}{\delta^3}z^4 + \frac{48}{\delta^2}z^3 - \frac{24}{\delta}z^2 + 5z - \delta, & \frac{\delta}{4} \leq z \leq \frac{\delta}{2}, \\ z - \delta, & z \geq \frac{\delta}{2}. \end{cases} \quad (3.7)$$

Then we have

$$\eta'(z) = \begin{cases} 0, & z \le \frac{\delta}{4}, \\ \frac{16}{\delta^2}(z - \frac{\delta}{4})^2(5 - \frac{8}{\delta}z), & \frac{\delta}{4} \le z \le \frac{\delta}{2}, \\ 1, & z \ge \frac{\delta}{2}. \end{cases} \quad (3.8)$$

$$\eta''(z) = \begin{cases} 0, & z \le \frac{\delta}{4}, \\ \frac{192}{\delta^2}(z - \frac{\delta}{4})(1 - \frac{2}{\delta}z), & \frac{\delta}{4} \le z \le \frac{\delta}{2}, \\ 0, & z \ge \frac{\delta}{2}. \end{cases} \quad (3.9)$$

$$\eta'''(z) = \begin{cases} 0, & z \le \frac{\delta}{4}, \\ \frac{96}{\delta^2}(3 - \frac{8}{\delta}z), & \frac{\delta}{4} \le z \le \frac{\delta}{2}, \\ 0, & z \ge \frac{\delta}{2}. \end{cases} \quad (3.10)$$

So

$$\eta \in C^2(R), \ 0 \le \eta'(z) \le 1 \ \text{and} \ 0 \le \eta''(z) \le \frac{6}{\delta}. \quad (3.11)$$

In this paper we only consider the $\eta(z)$ in (3.7) without loss of generality.

**Lemma 3.3** *Fix $0 < \delta < 2\sqrt{6t^*}$ ($t^*$ see (3.3)) then*

$$\eta'(z) = 0, \quad for \ z \le \tfrac{\delta}{4};$$

$$0 < \eta'(z) \le a_\delta < 1, \quad for \ z \in \left(\tfrac{\delta}{4}, \tfrac{3\delta}{8} + \tfrac{\delta}{8}\sqrt{1 - \tfrac{\delta^2}{24t^*}}\right]; \quad (3.12)$$

$$0 < \eta'(z) \le 1, \quad for \ z \in \left[\tfrac{3\delta}{8} + \tfrac{\delta}{8}\sqrt{1 - \tfrac{\delta^2}{24t^*}}, \tfrac{\delta}{2}\right];$$

$$\eta'(z) = 1, \quad for \ z \ge \tfrac{\delta}{2}.$$

**Lemma 3.4** *Let $0 < \delta < 2\sqrt{6t^*}$, assume that $\phi \in C^2(R^n \times [0, \infty))$ and $(x_0, t_0) \in (R^n \times (0, t^*))$ satisfy the following condition*

$$\eta(d) - \phi \quad (3.13)$$

*has a local minimum at the point $(x_0, t_0)$, then*

$$|\nabla \phi| = |\eta'(d)|, \quad at \ the \ point \ (x_0, t_0).$$

**Lemma 3.5** *Let $0 < \delta < 2\sqrt{6t^*}$ be fixed, and assume that $\phi \in C^2(R^n \times [0, \infty))$ and $(x_0, t_0) \in R^n \times (0, t^*)$ satisfy the following condition*

$$\eta(d) - \phi \ has \ a \ local \ minimum \ at \ the \ point \ (x_0, t_0), \quad (3.14)$$

*then*

$$\phi_t - \triangle\phi \ge 0, \quad for \ d(x_0, t_0) < \frac{\delta}{4}; \quad (3.15)$$

$$\phi_t - \triangle\phi \ge -\frac{6}{\delta}, \quad for \ d(x_0, t_0) \in \left[\frac{\delta}{4}, \frac{3\delta}{8} + \frac{\delta}{8}\sqrt{1 - \frac{\delta^2}{24t^*}}\right]; \quad (3.16)$$

$$\phi_t - \triangle \phi \geq -\frac{\delta}{4t^*}, \quad for \quad d(x_0, t_0) \in \left[\frac{3\delta}{8} + \frac{\delta}{8}\sqrt{1 - \frac{\delta^2}{24t^*}}, \frac{\delta}{2}\right]; \tag{3.17}$$

$$\phi_t - \triangle \phi \geq 0, \quad for \quad d(x_0, t_0) > \frac{\delta}{2}. \tag{3.18}$$

**Remark 3.6** *From Lemma 3.5, we have proved the result below in the sense of viscosity solution in [8]*

$$\eta(d)_t - \triangle \eta(d) \geq -\frac{C}{\delta}.$$

*where $C = 6$. Notice that the constructing of $\eta$ is more and more difficult along with the $C$ diminishing. We have to consider the very finer construction of $\eta$. This is the quite key for compatibility among those inequalities in our definitions.*

To find the solution of (1.1), it is sufficient to solve the following equations

$$\begin{cases} v_2 v_{1r} - v_1 v_{2r} = 0 \\ v_{3r}(r_t - \triangle r) - (v_{3rr} + |v_r|^2 v_3)|\triangledown r|^2 + \frac{1}{\varepsilon^2}v_3(1 - v_3^2)(\alpha + \beta - 2\alpha v_3^2) = 0. \end{cases}$$

But in order to find out a $\delta$-viscosity supersolution of (1.1), which are neither viscosity supersolution nor subsolution of (1.1), we need some lemmas globally control to $v$.

## 4 Global control maps

In order to control all solutions of (1.1) globally into the unit sphere, we first study the following ordinary differential equations:

$$\begin{cases} v_2(r)v_1''(r) - v_1(r)v_2''(r) = 0, \\ v_3''(r) + |v'(r)|^2 v_3(r) - \frac{\alpha+\beta-2\alpha v_3^2(r)}{\varepsilon^2} v_3(r)(1 - v_3^2(r)) = 0, \quad on\ R \\ v \in S^2, \end{cases} \tag{4.1}$$

**Lemma 4.1 (Global control map I)** *(1.1) possesses an explicit global control map from $R^{n+1}$ to $S^2 : v = v^\varepsilon(r) \in C^2(R; R^3) \cap S^2$ satisfying the ordinary differential equations (4.1)*

$$\begin{cases} v_1^\varepsilon(r) = \frac{(\alpha+\beta)\sin\varepsilon^2\cos(\varepsilon^2+\mu r\sin\varepsilon^2)-\sin(\mu r\sin\varepsilon^2)(\exp(2\mu r\cos\varepsilon^2)-0.5\beta)}{\sqrt{(\exp(2\mu r\cos\varepsilon^2)+\alpha+0.5\beta)^2-\alpha(\alpha+\beta)\cos\varepsilon^2}}, \\ v_2^\varepsilon(r) = \frac{(\alpha+\beta)\sin\varepsilon^2\sin(\varepsilon^2+\mu r\sin\varepsilon^2)+\cos(\mu r sin\varepsilon^2)(\exp(2\mu r\cos\varepsilon^2)-0.5\beta)}{\sqrt{(\exp(2\mu r\cos\varepsilon^2)+\alpha+0.5\beta)^2-\alpha(\alpha+\beta)\cos\varepsilon^2}}, \\ v_3^\varepsilon(r) = \frac{\sqrt{2(\alpha+\beta)\cos\varepsilon^2}\exp(\mu r\cos\varepsilon^2)}{\sqrt{(\exp(2\mu r\cos\varepsilon^2)+\alpha+0.5\beta)^2-\alpha(\alpha+\beta)\cos\varepsilon^2}}, \end{cases} \tag{4.2}$$

*where $\mu = \frac{\sqrt{\alpha+\beta}}{\varepsilon}$. Further $v^\varepsilon$ satisfies*

$$\begin{cases} \lim_{\varepsilon \to 0}\{r \in R \mid v_{3r}^\varepsilon(r) = 0\} = \{0\}, \\ |v_3^\varepsilon(r)| \neq 1\ for\ r \in R, \\ v_2^\varepsilon(r)(v_1^\varepsilon(r))' - v_1^\varepsilon(r)(v_2^\varepsilon(r))' = -\mu\sin\varepsilon^2 \neq 0\ for\ r \in R, \\ v_2^\varepsilon(r)(v_1^\varepsilon(r))'' - v_1^\varepsilon(r)(v_2^\varepsilon(r))'' = 0\ for\ r \in R, \end{cases} \tag{4.3}$$

as well as

$$\lim_{\varepsilon \to 0^+} v^\varepsilon(r) = \begin{cases} (0, -1, 0), & r > 0, \\ (0, \frac{2-\beta}{\sqrt{8\alpha+(2+\beta)^2}}, \frac{2\sqrt{2}\sqrt{\alpha+\beta}}{\sqrt{8\alpha+(2+\beta)^2}}), & r = 0, \\ (0, 1, 0), & r < 0. \end{cases} \quad (4.4)$$

Notice that we require $v^\varepsilon$ satisfying (4.3) and (4.4) for getting the existence and limiting behavior of $\delta$-viscosity supersolution of (1.1), respectively. Such requirements are most technical (see the proof of Theorem 4.1 in the next section).

Under the control map I (4.2), the equations (1.1) or (2.3) becomes

$$\begin{cases} -\frac{\sqrt{\alpha+\beta}}{\varepsilon} \sin \varepsilon^2 \, (r_t - \triangle r) = 0, \\ (v_3^\varepsilon(r))'(r_t - \triangle r) + \frac{1}{\varepsilon^2}(1 - |\nabla r|^2)v_3^\varepsilon(r)(1 - (v_3^\varepsilon(r))^2)(\alpha + \beta - 2\alpha v_3^2(r)) = 0, \, in \, \Omega \times (0, \infty). \end{cases} \quad (4.5)$$

This equation has only static solutions $r = c_0 + c_1 x_1 + c_2 x_2 + \cdots + c_n x_n$, where $c_1^2 + c_2^2 + \cdots + c_n^2 = 1$. This seems to show that the global control map I (4.2) is useless, but later we will see some global control maps, $\delta$-viscosity supersolutions (subsolutions), viscosity solutions and exact solutions constrained by an appropriate $\delta$-viscosity supersolution (subsolution) produced from the global control map I (4.2).

**Lemma 4.2 (Global control map II)** *Let $\beta > 0$, then (1.1) possesses an explicit global control map from $R^{n+1}$ to $S^2$ : $v = v^\varepsilon(r) \in C^2(R; R^3) \cap S^2$*

$$\begin{cases} v_1^\varepsilon(r) = \frac{\exp(2\mu r) - 0.5\beta}{\sqrt{(\exp(2\mu r) + 0.5\beta)^2 + 2\alpha \exp(2\mu r)}} \cos k, \\ v_2^\varepsilon(r) = \frac{\exp(2\mu r) - 0.5\beta}{\sqrt{(\exp(2\mu r) + 0.5\beta)^2 + 2\alpha \exp(2\mu r)}} \sin k, \\ v_3^\varepsilon(r) = \frac{\sqrt{2(\alpha+\beta)} \exp(\mu r)}{\sqrt{(\exp(2\mu r) + 0.5\beta)^2 + 2\alpha \exp(2\mu r)}}, \end{cases} \quad (4.6)$$

*where $\mu = \frac{\sqrt{\alpha+\beta}}{\varepsilon}$, $k$ is any given constant. Further $v^\varepsilon$ satisfies*

$$\begin{cases} \lim_{\varepsilon \to 0}\{r \in R \mid v_{3r}^\varepsilon(r) = 0\} = \{0\}, \\ |v_3^\varepsilon(r)| \neq 1, & \text{for } r \in R, \\ v_2^\varepsilon(r)(v_1^\varepsilon(r))' - v_1^\varepsilon(r)(v_2^\varepsilon(r))' = 0, & \text{for } r \in R, \end{cases} \quad (4.7)$$

as well as

$$\lim_{\varepsilon \to 0^+} v^\varepsilon(r) = \begin{cases} (\cos k, \sin k, 0), & r > 0, \\ (\frac{1-0.5\beta}{\sqrt{(1+0.5\beta)^2+2\alpha}} \cos k, \frac{1-0.5\beta}{\sqrt{(1+0.5\beta)^2+2\alpha}} \sin k, \frac{\sqrt{2(\alpha+\beta)}}{\sqrt{(1+0.5\beta)^2+2\alpha}}), & r = 0, \\ (-\cos k, -\sin k, 0), & r < 0. \end{cases} \quad (4.8)$$

Under the global control map II (4.6) $u = v^\varepsilon(r)$, $r = r(x,t)$, the first equation of (2.3) becomes an identical relation, and (2.3) becomes just the second equation of the equation (4.5). Thus (1.1) becomes

$$(v_3^\varepsilon(r))'(r_t - \triangle r) + \frac{1}{\varepsilon^2}(1 - |\nabla r|^2)v_3^\varepsilon(r)(1 - (v_3^\varepsilon(r))^2)(\alpha + \beta - 2\alpha v_3^2(r)) = 0, \quad in \ \Omega \times (0, \infty). \quad (4.9)$$

We find that the asymptotic behavior of the solution under global control map II is very similar to the one under global control map I. These solutions are constrained by a $\delta$-viscosity supersolution and subsolution. As $\delta$-viscosity supersolution and subsolution are approximate solution, them can often be explicitly given so that our study feasible. When we consider the approximate solution of (1.1) in $C^{2+1}(\bar{\Omega} \times [0,T])$ under global control map I, we only need to solve the equation (4.9) as $\varepsilon$ small enough.

## 5 Construction of $\delta$-viscosity supersolution and subsolution

In this section we apply the signed distance function $d$ to build $\delta$-viscosity solution of (1.1), they are neither viscosity supersolution nor viscosity subsolution of (1.1). Now and afterwards we merely discuss $\delta$-viscosity supersolution, and the argument is similar for $\delta$-viscosity subsolution. Now we state and prove one of our prime results, theorem on $\delta$-viscosity supersolution of (1.1).

**Theorem 5.1** *For every fixed $0 < \delta < 2\sqrt{6t^*}$, there exists $\varepsilon_0 = \varepsilon_0(\delta) > 0$ such that $u = v^\varepsilon(r)$ is a $\delta$-viscosity supersolution (respectively, subsolution) of (1.1) on $R^n \times (0, t^*)$ with initial value $r(x,0) = \eta(d(x,0))$ (respectively, $-\eta(-d(x,0))$), where $r = \eta(d(x,t)) + \frac{\delta}{4t^*}t$ (respectively, $-\eta(-d(x,t)) - \frac{\delta}{4t^*}t$) and $0 < \varepsilon \leq \varepsilon_0$.*

**Remark 5.2** $\phi_t - \triangle\phi < 0$, when $d(x_0,t_0) \in \left[\frac{\delta}{4}, \frac{3\delta}{8} + \frac{\delta}{8}\sqrt{1 - \frac{\delta^2}{24t^*}}\right]$. $u = v^\varepsilon(\eta(d(x,t)) + \frac{\delta}{4t^*}t)$ is not a supersolution

**Remark 5.3** *The inequalities (3.16) and (3.17) in the proof of Lemma 3.5 cannot be improved to the following inequality:*
$$\phi_t - \triangle\phi \geq 0, \qquad d(x_0, t_0) \in (\tfrac{\delta}{4}, \tfrac{\delta}{2}).$$
*In fact the $v^\varepsilon(r)$ in Theorem 5.1 is not viscosity supersolution of the first equation of (2.3). Hence the $v^\varepsilon(r)$ is not viscosity supersolution of (1.1). Clearly $v^\varepsilon(r)$ is not also viscosity subsolution of (1.1).*

## 6 Maximum and comparison principle under the global control map, existence and uniqueness of $\delta$-viscosity solution

**Theorem 6.1 (Maximum principle under the global control map I, II from $R^{n+1}$ to $S^2$)** Let $0 < \delta < 2\sqrt{6t^*}$ be fixed, $U_M = \{x \in R^n |\ |x| < M\}$.

(i) *If three dimensional vector valued function $u = v^\varepsilon(r) \in C^2(R; R^3) \cap S^2$ is a $\delta$-viscosity supersolution of (1.1). Assume that there is a constant $M$ and a point $(\widetilde{x}, \widetilde{t}) \in U_M \times [0, t^*]$ such that $r(\widetilde{x}, \widetilde{t}) \leq r(x, t)$ for every pair $x \notin U_M$ and $t \in [0, t^*]$. Then*

$$\min_{(x,t) \in R^n \times [0,t^*]} r(x,t) = \min_{x \in \overline{U}_M} r(x,0), \quad or \quad \frac{\delta}{4} \leq r(x,t) \leq \frac{\delta}{2}.$$

(ii) *If three dimensional vector valued function $u = v^\varepsilon(r) \in C^2(R; R^3) \cap S^2$ is a $\delta$-viscosity subsolution of (1.1). Assume that there is a constant $M$ and a point $(\widetilde{x}, \widetilde{t}) \in U_M \times [0, t^*]$ such that*

$r(\widetilde{x},\widetilde{t}) \geq r(x,t)$ for every pair $x \notin U_M$ and $t \in [0,t^*]$. Then

$$\max_{(x,t)\in R^n \times [0,t^*]} r(x,t) = \max_{x \in \overline{U}_M} r(x,0), \quad or \quad -\frac{\delta}{2} \leq r(x,t) \leq -\frac{\delta}{4}.$$

Applying the above maximum principle we obtain the following comparison principle.

**Theorem 6.2 (Comparison principle)** *Let $0 < \delta < 2\sqrt{6t^*}$ be fixed, $U_M = \{x \in R^n| |x| < M\}$. If three dimensional vector valued function $v^\varepsilon(r) \in C^2(R;R^3) \cap S^2$ is a $\delta$-viscosity supersolution and $v^\varepsilon(s) \in C^2(R;R^3) \cap S^2$ is a $\delta$-viscosity subsolution of (1.1). Assume that there is a constant $M > 0$ and $(\widetilde{x},\widetilde{t}) \in U_M \times [0,t^*]$ such that $s(x,t) \leq r(\widetilde{x},\widetilde{t}) \leq r(x,t)$ for every pair $x \notin U_M$ and $t \in \times[0,t^*]$. Suppose further*

$$s \leq r \ on \ R^n \times \{t=0\}.$$

*If $r$ and $s$ continue, then*

$$s \leq r \quad on \ R^n \times [0,t^*] \quad when \ r \notin [\frac{\delta}{4},\frac{\delta}{2}] \ and \ s \notin [-\frac{\delta}{2},-\frac{\delta}{4}].$$

Since $\delta$-viscosity solution is a viscosity of one equation of being considered system, similar to the argument of [9] to the existence and uniqueness of the viscosity solution, we draw the following existence and uniqueness theorem.

**Theorem 6.3** *Let $0 < \delta < 2\sqrt{6t^*}$ be fixed. Assume that there exist some $M > 0$ such that*

$$r_0(x) = |x| \ for \ every \ |x| \geq M.$$

*If $r_0(x) \in C(R^n)$, then under the global control map I or II there is a unique $r(x,t) \in C(R^n \times [0,t^*])$ such that three dimensional vector valued function $v^\varepsilon(r) \in C^2(R;R^3) \cap S^2$ is a $\delta$-viscosity solution of the equations (1.1) with initial value $r(x,0) = r_0(x)$.*

On the one hand the proof on this theorem is similar to the proof on the viscosity solution in [9], on the other hand the asymptotic behavior considered by us can be controlled only by a pair $\delta$-viscosity supersolution and $\delta$-viscosity subsolution, so we omit the proof.

Let us consider the following initial problem from (4.9).

$$\begin{cases} (v_3^\varepsilon(r))'(r_t - \triangle r) + \frac{1}{\varepsilon^2}(1 - |\nabla r|^2)v_3^\varepsilon(r)(1 - (v_3^\varepsilon(r))^2)(\alpha + \beta - 2\alpha v_3^2(r)) = 0, \ in \ R^n \times (0,\infty), \\ r(x,0) = r_0 \in C(R^n). \end{cases} \tag{6.1}$$

(6.1) farther becomes

$$\begin{cases} r_t - \triangle r = 2\mu(\alpha+\beta)(1 - |\nabla r|^2)\frac{e^{\mu r}-0.5\beta e^{-\mu r}}{e^{\mu r}+0.5\beta e^{-\mu r}}\frac{(e^{\mu r}+0.5\beta e^{-\mu r})^2 - 2\alpha}{(e^{\mu r}+0.5\beta e^{-\mu r})^2 + 2\alpha}, \\ r(x,0) = r_0. \end{cases} \tag{6.2}$$

Combining Theorem 6.1, Theorem 6.2 and Theorem 6.3 we have the comparison principle.

**Theorem 6.4** *Let $0 < \delta < 2\sqrt{6t^*}$ be fixed. Assume $r_0 \in C^{2,\alpha}(R^n)$ and that there exist some $M > 0$ such that*

$$r_0(x) = |x| \ for \ every \ |x| \geq M.$$

*Then the solution of (6.2)*

$$r \in C^{2,\alpha;\,1,\frac{\alpha}{2}}(R^n \times [0,t^*]) \cap C^\infty(R^n \times (0,t_0))$$

*can be controlled a pair of $\delta$-viscosity supersolution $r^+$ and $\delta$-viscosity subsolution $r^-$:*

$$r^- \leq r \leq r^+, \text{ for } (x,t) \in R^n \times [0,t_*], t_* = \min\{t^*, t_0\}.$$

# 7 Asymptotic analysis

Let $d_0$ be the signed distance function to $\Gamma_0$, and set

$$g(x) = d_0(x), \ x \in R^n. \tag{7.1}$$

In this section we consider the asymptotic behavior of the $\delta$-viscosity supersolution and subsolution of (1.1) with initial value

$$u(x,0) = v^\varepsilon(\eta(g)).$$

We however always reset $t_* = \min\{t^*, t_0\}$, $r_+ = \eta(d(x,t)) + \frac{\delta}{4t^*}t$, $r_- = -\eta(-d(x,t)) - \frac{\delta}{4t^*}t$, $I_* \equiv \{(x,t) \in R^n \times (0,t_*) \mid w(x,t) > 0\}$, $O_* \equiv \{(x,t) \in R^n \times (0,t_*) \mid w(x,t) < 0\}$, the interface between $I_*$ and $O_*$: $\Gamma_* \equiv \{(x,t) \in R^n \times (0,t_*) \mid w(x,t) = 0\}$. The interface is a generalized motion governed by mean curvature. Here, $w$ is given in (3.2).

**Theorem 7.1** *Let $\delta > 0$ is small enough, then*

$$u = v^\varepsilon(r_\pm) \to (0,1,0) \text{ uniformly on every compact subset of } I_*, \tag{7.2}$$

$$u = v^\varepsilon(r_\pm) \to (0,-1,0) \text{ uniformly on every compact subset of } O_*, \tag{7.3}$$

*and*

$$u = v^\varepsilon(r_\pm) \to \vec{a} \text{ somewhere near } \Gamma_*, \tag{7.4}$$

*where $\vec{a} = (0, \frac{2-\beta}{\sqrt{8\alpha+(2+\beta)^2}}, \frac{2\sqrt{2}\sqrt{\alpha+\beta}}{\sqrt{8\alpha+(2+\beta)^2}})$.*

We have the following results:

**Theorem 7.2** *Let $0 < \delta < 2\sqrt{6t^*}$ be fixed. Assume that $r_0 \in C^{2,\alpha}(R^n)$ and there exist some $M > 0$ such that*

$$r_0(x) = |x| \text{ for every } |x| \geq M.$$

*Then for every fixed $\varepsilon > 0$, $\exists\, t_0 > 0$ such that (6.2) has a classical solution*

$$r = \omega \in C^{2,\alpha;\,1,\frac{\alpha}{2}}(R^n \times [0,t^*]) \cap C^\infty(R^n \times (0,t_0)).$$

*Furthermore (1.1) has a classical solution $u = v^\varepsilon(r) \in C^2(R^n \times [0,t_*]; R^3) \cap S^2$ satisfying for $A$ and $B$ be arbitrary given compact subsets in $I_*$ and in $O_*$, respectively, there exist $\delta > 0$ such that*

$$u = v^\varepsilon(r) \to \begin{cases} (0,1,0) \text{ uniformly on } A, \\ (0,-1,0) \text{ uniformly on } B \end{cases}$$

*as $\varepsilon \to \infty$, where $t_* = \min\{t^*, t_0\}$.*

**Remark 7.3** *Multiply $v^\varepsilon$ in (4.2) by an orthogonal constant matrix of order three with determinant value 1 we may construct other solution of (4.1). Consequently other $\delta$-viscosity supersolutions and subsolutions can be constructed.*

# 8 Direct method of viscosity solution

In Section 5, we constructed a class of $\delta$-viscosity supersolution of (1.1). Although it is not viscosity supersolution for (1.1), but using it we can show the existence of viscosity supersolution and classical solution of (1.1) for a class of initial conditions. The $\delta$-viscosity supersolution enlighten us to obtain an explicit viscosity supersolution of (1.1) with special initial value (re. Theorem 8.2, Remark 8.3).

**Theorem 8.1** $u = v^\varepsilon(\eta(d) + \frac{\delta}{4t^*}t)$ is a viscosity supersolution of (1.1).

**Theorem 8.2** If $-\frac{1}{4}\beta < \alpha < \frac{\sqrt{7}-1}{4}\beta$ and $\beta > 0$, then for every constant $k$, the equations (1.1) with initial

$$u_0^\varepsilon = (v_1^\varepsilon(\eta(d_0)), v_2^\varepsilon(\eta(d_0)), v_3^\varepsilon(\eta(d_0))) \tag{8.1}$$

has a classical solution $u^\varepsilon$ satisfying

$$u^\varepsilon \to (\cos k, \sin k, 0) \; uniformly \; on \; compact \; subsets \; of \; I_*, \tag{8.2}$$

$$u^\varepsilon \to (-\cos k, -\sin k, 0) \; uniformly \; on \; compact \; subsets \; of \; O_*. \tag{8.3}$$

**Remark 8.3** It is known that if $\Gamma_0$ is the sphere $\partial B(0, R)$, then

$$\Gamma_t = \begin{cases} \partial B(0, R(t)), & if \; 0 \leq t < t^*, \\ 0, & if \; t = t^*, \\ \phi, & if \; t > t^*, \end{cases} \tag{8.4}$$

where $R(t) = \sqrt{R^2 - 2(n-1)t}$ for $0 \leq t \leq t^* = \frac{R^2}{2(n-1)}$. If $v^\varepsilon$ is taken as the one defined in Lemma 4.1( 4.2, respectively), then $u = v^\varepsilon(\eta(d(x,t)) + \frac{\delta}{4t^*}t)$ is a $\delta$-viscosity (viscosity, respectively) supersolution of the first equation of (1.1) by Theorem 5.1(8.1, respectively), where $d(x,t) = R(t) - \sqrt{x_1^2 + x_2^2 + \cdots + x_n^2}$.

**Remark 8.4** The second equation of (2.1) differs from Allen-Cahn equation, because it submit to the first equation of (2.1), $u \in S^2$ and the gradient term. But if let $\omega = \sin^{-1}(v_3^\varepsilon(r))$ then (1.1) becomes the Allen-Cahn equation with the positive kinetic constant $\alpha = \frac{1}{\varepsilon^2}$, gradient energy coefficient $k = 0.5\varepsilon^2$ and the nonlinear term $\sin\omega\cos\omega(\alpha + \beta - 2\alpha\sin^2\omega)$:

$$\omega_t - \Delta\omega + \tfrac{1}{\varepsilon^2} \sin\omega\cos\omega(\alpha + \beta - 2\alpha\sin^2\omega) = 0. \tag{8.5}$$

Notice that when we consider static solution, $u_3 = \sin\omega = v_3^\varepsilon(r)$ do not satisfy the condition "g can be smoothly extended to a function from $\overline{\Omega}$ to $S^2 \cap \{u_3 = 0\}$" in Theorem A. If $\alpha = 0$ then $v_3^\varepsilon(r)$ is strictly increasing, and there exists its inverse function. By Theorem 8.2, the equation (8.5) has a solution satisfying (8.1) and (8.2). If we take $k = 1$ in Allen-Cahn equation, then only when $\varepsilon = \sqrt{2}$ it correspond to (8.5).

# 9 The appendix on Lemma4.1

We verify
$$v_2 v_1' - v_1 v_2' = -\mu \sin \epsilon^2;$$
$$v_3'' + \frac{\mu^2 \sin^2 \epsilon^2 + (v_3')^2}{1-v_3^2} v_3 - \frac{1}{\varepsilon^2} v_3 (1-v_3^2)(\alpha + \beta - 2\alpha v_3^2) = 0;$$
$$\lim_{\varepsilon \to 0} \{r \in R \mid v_{3r}(r) = 0\} = \{0\}.$$

Let $A = \sqrt{\alpha + \beta} \cos(\epsilon^2)$, $B = \sqrt{\alpha + \beta} \sin(\epsilon^2)$.

$$v_{1r} = \frac{-2 e^{\frac{2rA}{\epsilon}} \left(e^{\frac{2rA}{\epsilon}} + \alpha + \frac{\beta}{2}\right) A \left((\alpha+\beta) \cos(\epsilon^2 + \frac{rB}{\epsilon}) \sin(\epsilon^2) - \left(e^{\frac{2rA}{\epsilon}} - \frac{\beta}{2}\right) \sin(\frac{rB}{\epsilon})\right)}{\epsilon \left(\left(e^{\frac{2rA}{\epsilon}} + \alpha + \frac{\beta}{2}\right)^2 - \alpha(\alpha+\beta) \cos(\epsilon^2)^2\right)^{\frac{3}{2}}}$$
$$- \frac{\sqrt{\alpha+\beta} \left(\left(2 e^{\frac{2rA}{\epsilon}} - \beta\right) \cos(\frac{rB}{\epsilon}) \sin(\epsilon^2) + 2 \left(2 e^{\frac{2rA}{\epsilon}} \cos(\epsilon^2) \sin(\frac{rB}{\epsilon}) + (\alpha+\beta) \sin(\epsilon^2)^2 \sin(\epsilon^2 + \frac{rB}{\epsilon})\right)\right)}{2\epsilon \sqrt{\left(e^{\frac{2rA}{\epsilon}} + \alpha + \frac{\beta}{2}\right)^2 - \alpha(\alpha+\beta) \cos(\epsilon^2)^2}},$$

$$v_{2r} = \frac{\sqrt{\alpha+\beta} \left(4 e^{\frac{2rA}{\epsilon}} \cos(\epsilon^2) \cos(\frac{rB}{\epsilon}) + \sin(\epsilon^2) \left(2(\alpha+\beta) \cos(\epsilon^2 + \frac{rB}{\epsilon}) \sin(\epsilon^2) - \left(2 e^{\frac{2rA}{\epsilon}} - \beta\right) \sin(\frac{rB}{\epsilon})\right)\right)}{2\epsilon \sqrt{\left(e^{\frac{2rA}{\epsilon}} + \alpha + \frac{\beta}{2}\right)^2 - \alpha(\alpha+\beta) \cos(\epsilon^2)^2}}$$
$$- \frac{2 e^{\frac{2rA}{\epsilon}} \left(e^{\frac{2rA}{\epsilon}} + \alpha + \frac{\beta}{2}\right) A \left(\left(e^{\frac{2rA}{\epsilon}} - \frac{\beta}{2}\right) \cos(\frac{rB}{\epsilon}) + (\alpha+\beta) \sin(\epsilon^2) \sin(\epsilon^2 + \frac{rB}{\epsilon})\right)}{\epsilon \left(\left(e^{\frac{2rA}{\epsilon}} + \alpha + \frac{\beta}{2}\right)^2 - \alpha(\alpha+\beta) \cos(\epsilon^2)^2\right)^{\frac{3}{2}}},$$

$$v_{3r} = \frac{-2\sqrt{2} e^{\frac{rA}{\epsilon}} (\alpha+\beta) \cos(\epsilon^2)^2 \left(4 e^{\frac{4rA}{\epsilon}} - 2\alpha^2 - 2\alpha\beta - \beta^2 + 2\alpha(\alpha+\beta) \cos(2\epsilon^2)\right)}{\epsilon \left(4 e^{\frac{4rA}{\epsilon}} + 2\alpha^2 + 2\alpha\beta + \beta^2 + 4 e^{\frac{2rA}{\epsilon}} (2\alpha+\beta) - 2\alpha(\alpha+\beta) \cos(2\epsilon^2)\right)^{\frac{3}{2}}},$$

$$v_{3rr} = \frac{-4\sqrt{2} e^{\frac{3rA}{\epsilon}} \left(e^{\frac{2rA}{\epsilon}} + \alpha + \frac{\beta}{2}\right) (\alpha+\beta)^{\frac{3}{2}} \cos(\epsilon^2)^3}{\epsilon^2 \left(\left(e^{\frac{2rA}{\epsilon}} + \alpha + \frac{\beta}{2}\right)^2 - \alpha(\alpha+\beta) \cos(\epsilon^2)^2\right)^{\frac{3}{2}}} + \frac{\sqrt{2} e^{\frac{rA}{\epsilon}} (\alpha+\beta)^{\frac{3}{2}} \cos(\epsilon^2)^3}{\epsilon^2 \sqrt{\left(e^{\frac{2rA}{\epsilon}} + \alpha + \frac{\beta}{2}\right)^2 - \alpha(\alpha+\beta) \cos(\epsilon^2)^2}}$$
$$+ \frac{16\sqrt{2} e^{\frac{3rA}{\epsilon}} (\alpha+\beta)^{\frac{3}{2}} \cos(\epsilon^2)^3 \left(8 e^{\frac{6rA}{\epsilon}} - 4\alpha^3 - 6\alpha^2\beta - 2 e^{\frac{2rA}{\epsilon}} \beta^2 - 4\alpha\beta^2 - \beta^3 + 4 e^{\frac{4rA}{\epsilon}} (2\alpha+\beta) + 2\alpha(\alpha+\beta) \left(4 e^{\frac{2rA}{\epsilon}} + 2\alpha + \beta\right) \cos(2\epsilon^2)\right)}{\epsilon^2 \left(4 e^{\frac{4rA}{\epsilon}} + 2\alpha^2 + 2\alpha\beta + \beta^2 + 4 e^{\frac{2rA}{\epsilon}} (2\alpha+\beta) - 2\alpha(\alpha+\beta) \cos(2\epsilon^2)\right)^{\frac{5}{2}}}$$

So $v_2 v_1' - v_1 v_2' = -\frac{\sqrt{\alpha+\beta} \sin \epsilon^2}{\epsilon}$.

$\frac{(-\frac{\sqrt{\alpha+\beta} \sin \epsilon^2}{\epsilon})^2 + v_{3r}^2}{1-v_3^2}$

$$= \frac{\sqrt{2} e^{\frac{rA}{\epsilon}} (\alpha+\beta)^{\frac{3}{2}} \cos(\epsilon^2) \left(\frac{8 e^{\frac{2rA}{\epsilon}} (\alpha+\beta) \cos(\epsilon^2)^4 \left(-4 e^{\frac{4rA}{\epsilon}} + 2\alpha^2 + 2\alpha\beta + \beta^2 - 2\alpha(\alpha+\beta) \cos(2\epsilon^2)\right)^2}{\left(4 e^{\frac{4rA}{\epsilon}} + 2\alpha^2 + 2\alpha\beta + \beta^2 + 4 e^{\frac{2rA}{\epsilon}} (2\alpha+\beta) - 2\alpha(\alpha+\beta) \cos(2\epsilon^2)\right)^3} + \sin(\epsilon^2)^2\right)}{\epsilon^2 \sqrt{\left(e^{\frac{2rA}{\epsilon}} + \alpha + \frac{\beta}{2}\right)^2 - \alpha(\alpha+\beta) \cos(\epsilon^2)^2} \left(1 - \frac{8 e^{\frac{2rA}{\epsilon}} (\alpha+\beta) \cos(\epsilon^2)^2}{\left(2 e^{\frac{2rA}{\epsilon}} + 2\alpha+\beta\right)^2 - 4\alpha(\alpha+\beta) \cos(\epsilon^2)^2}\right)},$$

$v_3(1-v_3^2)(\alpha+\beta - 2\alpha v_3^2)$

$$= \frac{\sqrt{2} e^{\frac{rA}{\epsilon}} \sqrt{\alpha+\beta} \cos(\epsilon^2) \left(1 - \frac{8 e^{\frac{2rA}{\epsilon}} (\alpha+\beta) \cos(\epsilon^2)^2}{\left(2 e^{\frac{2rA}{\epsilon}} + 2\alpha+\beta\right)^2 - 4\alpha(\alpha+\beta) \cos(\epsilon^2)^2}\right) \left(\alpha+\beta - \frac{4 e^{\frac{2rA}{\epsilon}} \alpha(\alpha+\beta) \cos(\epsilon^2)^2}{\left(e^{\frac{2rA}{\epsilon}} + \alpha + \frac{\beta}{2}\right)^2 - \alpha(\alpha+\beta) \cos(\epsilon^2)^2}\right)}{\sqrt{\left(e^{\frac{2rA}{\epsilon}} + \alpha + \frac{\beta}{2}\right)^2 - \alpha(\alpha+\beta) \cos(\epsilon^2)^2}}.$$

So
$$v_3'' + \frac{\mu^2 \sin^2 \epsilon^2 + (v_3')^2}{1-v_3^2} v_3 - \frac{1}{\varepsilon^2} v_3(1-v_3^2)(\alpha+\beta - 2\alpha v_3^2) = 0.$$

From the expression of $v_{3r}$, we have
$$\{r \in R \mid v_{3r}(r) = 0\} = \{\frac{\ln[(\alpha+0.5\beta)^2 \sin^2 \varepsilon^2 + 0.25\beta^2 \cos^2 \varepsilon^2]}{4\sqrt{\alpha+\beta} \cos \varepsilon^2} \varepsilon\}.$$

Thus
$$\lim_{\varepsilon \to 0}\{r \in R \mid v_{3r}(r) = 0\} = \{0\}.$$

**Acknowledgements** The first author is grateful to Professor Zhouping Xin of Chinese University of Hongkong, Professor Baoxiang Wang of Peking University and Professor Xiangao Liu of Fudan University for all of their guidance, support, and encouragement.

**Summary of the first author** The first author mainly studies nonlinear analysis and nonlinear partial differential equation. He visited the Institute of Mathematical Sciences of Chinese University of Hongkong, Institute of Mathematics of Fudan University and Institute of Mathematics of Peking University; published some papers in Math. Numer. Sinica(1988), Appl. Math. Mech.(1990), Proc. Amer. Math. Soc.(1998), Commun. Nonlinear. Sci. numer. simul(2000), J. Math. Phy.(2001), Physics Letters A(2003), Acta Mathematica Sinica(2004), J. Partial Diff. Eqs.(2007), Chinese Phys. Lett(2008), Acta Mathematicae Applicatae Sinica(2009), Nonlinear Analysis(2009), Dynamics of Continuous, Discrete and Impulsive Systems Series A(2010), Science China Mathematics(2011), Physics Letters A(2012), etc.